\documentclass[12pt,twoside]{article}

\usepackage{amsfonts}
\usepackage{amssymb}
\usepackage{latexsym}
\usepackage{amsbsy}
\usepackage{fullpage}
\usepackage[thmmarks]{ntheorem}
\usepackage{ifthen}

\theoremstyle{plain}
\theorembodyfont{\slshape}
\theoremseparator{.}
\newtheorem{THM}{\indent Theorem}[section]
\newtheorem{LEM}{\indent Lemma}[section]
\theorembodyfont{\normalfont}
\theoremstyle{nonumberplain}
\newtheorem{COR}{\indent Corollary}
\newtheorem{RMK}{\indent Remark}

\theoremsymbol{$\Box$}
\newtheorem{PRF}{\indent Proof}

\begin{document}

\parindent=1em
\parskip=0ex
\raggedbottom

\setlength{\headsep}{5ex}

\pagestyle{myheadings}
\thispagestyle{empty}

\makeatletter
\renewcommand{\section}{\@startsection
    {section}%
    {1}%
    {0em}%
    {2\baselineskip}%
    {\baselineskip}%
    {\normalfont\normalsize\bf}}%
\makeatother

\setcounter{section}{0}

\makeatletter
\renewcommand{\@seccntformat}[1]{\csname the#1\endcsname.\hspace{0.5em}}
\makeatother

\makeatletter
\renewcommand{\@cite}[2]{[{\bf #1}\ifthenelse{\boolean{@tempswa}}{, #2}{}]}
\makeatother

\makeatletter
\renewcommand{\@biblabel}[1]{#1.}
\makeatother

\markboth%
{\it M. J. Crabb, J. Duncan and C. M. McGregor \hfil}%
{\hfil \it On one-sided primitivity of Banach algebras}

\newcommand\suBset{\subseteq}
\newcommand{\LIM}{\mathop {\mathrm {LIM}} \limits}
\newcommand{\SPEC}{\mathop {\mathrm {Sp}} \nolimits}
\renewcommand{\epsilon}{\varepsilon}
\renewcommand{\theequation}{\thesection.\arabic{equation}}
\newcommand\SP{\hspace{0.2em}}
\newcommand\SQ{}
\newcommand\BS{\hspace{-0.001em}}

\begin{center}
{\large ON ONE-SIDED PRIMITIVITY OF BANACH ALGEBRAS} \\ [6ex] {M. J.
CRABB$^1$, J. DUNCAN$^2$ AND C. M. MCGREGOR$^1$} \\ [2ex]
$^1${\it Department of Mathematics, University of Glasgow,} \\
{\it Glasgow, G12 8QW, Scotland, UK} (m.crabb@gla.ac.uk, c.mcgregor@maths.gla.ac.uk) \\
$^2${\it Department of Mathematical Sciences, SCEN 301, University of Arkansas,} \\
{\it Fayetteville, AR 72701} (jduncan@uark.edu)
\end{center}

\vspace{6ex}
{\small\noindent
{\it Abstract}\quad
Let $S$ be the semigroup with identity, generated by $x$ and $y$,
subject to $y$ being invertible and $yx=xy^2$.
We study two Banach algebra completions of the semigroup algebra $\mathbb{C}S$.
Both completions are shown to be left-primitive and have separating families of
irreducible infinite-dimensional right modules.
As an appendix, we offer an alternative proof that $\mathbb{C}S$ is left-primitive
but not right-primitive.
We show further that, in contrast to the completions,
every irreducible right module for $\mathbb{C}S$ is finite dimensional
and hence that $\mathbb{C}S$ has a separating family of such modules.
}

\vspace{3ex}
{\small\noindent
{\it Keywords:}
Banach algebra;
primitive;
left-primitive;
right-primitive;
irreducible representation}

\vspace{3ex}
{\small\noindent
2000 {\it Mathematics subject classification:}
\begin{tabular}[t]{l}
% Primary 00A00; 00B00; 00C00 \\
% Secondary 00X00; 00Y00; 00Z00
Primary 46H20
\end{tabular}}

\section{Introduction}
\label{sec intro} \noindent In 1964 Bergman \cite{Be} gave the first
example of a ring which is right-primitive (has a faithful
irreducible right module) but not left-primitive. Further examples
were given in 1969 by Jategaonkar \cite{Ja}.  In 1979 Irving
\cite{Ir2} gave a whole family of examples of rings, even linear
algebras, which are left-primitive but not right-primitive.
The question of whether
or not there exists a Banach algebra with this property
was raised by Bonsall and Duncan~\cite{BoDu1}, and remains open.
With a view to examining the Banach
algebra situation we look in detail at one of the examples of
Irving, namely the semigroup algebra $\mathbb{C}S$ where $S$ is the
semigroup with identity, generated by $x$ and $y$, subject to $y$
being invertible and $yx=xy^2$. Thus
$$
S = \left\{ x^m y^n : m \in \mathbb{Z}_+,\ n \in \mathbb{Z},\ yx=xy^2 \right\},
$$
where $\mathbb{Z}_+ = \{k\in\mathbb{Z}:k\geq 0\}$.
In Section~\ref{sec appendix} we present a
proof of Irving's result for $\mathbb{C}S$ that may be more amenable
to Banach algebraists.   In fact, we prove even more; we prove that
every irreducible right module of $\mathbb{C}S$ is
finite-dimensional. Since $\mathbb{C}S$ is left-primitive and hence
semi-simple, it follows that $\mathbb{C}S$ has a separating family
of irreducible right representations on finite-dimensional spaces.
In the language of representation theory this says that $S$ is
residually finite. It follows that the algebra $\mathbb{C}S$ has the
property of direct finiteness, that is, for $a,b \in \mathbb{C}S$ we
have $ba=1$ if and only if $ab=1$.

Our principal aim in Section~\ref{sec rep thy B alg} is to study two Banach algebra
completions of $\mathbb{C}S$.   The obvious completion, $A$, say, is
given by $\ell^1(S)$ where we replace finite sums by absolutely
convergent sums with the usual $\ell^1$-norm.   It will be
convenient to view $A$ as a graded algebra; thus
$$
A = \Big\{a=\sum_{m=0}^\infty x^m \phi_m(y) :
\phi_m \in \ell^1(\mathbb{Z}),\ \sum_{m=0}^\infty \|\phi_m\|_1 < \infty \Big\}
$$
the norm being given by $\|a\| = \sum_{m=0}^\infty
\|\phi_m(y)\|_1$. The convolution Banach algebra
$\ell^1(\mathbb{Z})$ is isometrically isomorphic to the Wiener
algebra, $W$, of all continuous complex functions on $\mathbb{T}$,
the unit circle, with absolutely convergent Fourier series (the norm
being the absolute sum of the Fourier coefficients).  On occasions
it will be helpful to regard the $\phi_m$ as functions on $\mathbb{T}$.
We shall also consider the completion $B$ given by
$$
B = \Big\{ b=\sum_{m=0}^\infty x^m \phi_m :
\phi_m \in C(\mathbb{T}),\ \sum_{m=0}^\infty \|\phi_m\|_\infty < \infty \Big\}
$$
the norm
being given by $\|b\| = \sum_{m=0}^\infty \|\phi_m\|_\infty $. The
product in $A$ and in $B$ is determined by the formula
$$
x^m\phi(y) \SP x^n \psi(y) = x^{m+n} \phi(y^{2^n}) \psi(y).
$$
We shall show that both $A$ and $B$ are left-primitive and also residually
finite and hence also have the property of direct finiteness.  We
show that, in contrast to the situation for $\mathbb{C}S$, the
Banach algebras $A$ and $B$ have many non-faithful irreducible
infinite-dimensional right modules; enough, in fact, to separate the
points of the algebra. We are still unable to prove that either $A$
or $B$ fails to be right-primitive. A key step in the proof that
$\mathbb{C}S$ is not right primitive involves a technique that has
no analogue for {\it infinite\/} series.  It is amusing to note that
any answer to the right-primitivity of $A$ or $B$ will be
interesting.  If either $A$ or $B$ is {\it not\/} right-primitive, then we
have a desired example.  If they {\it are\/} right-primitive, then we have a
dramatic difference between the purely algebraic $\mathbb{C}S$ and
two of its natural Banach algebra completions.

We remark here that it is much easier to establish primitivity than
to establish non-primitivity.  To prove that an algebra $\mathcal{A}$ is
right-primitive it is enough to construct one faithful irreducible
right module.  To prove that $\mathcal{A}$ is not right-primitive we have to
show that {\it every\/} irreducible right module fails to be
faithful.  Equivalently, we have to identify every maximal modular
right ideal $K$ of $\mathcal{A}$ and show that the quotient ideal $K:\mathcal{A}$ is
always non-zero.  Bergman was able to do this for his example.
For Banach algebras, identifying all maximal modular right ideals
is usually a hopeless task.
We are thus forced (along with Irving)
to assume the existence of some faithful irreducible right
module and look for some contradiction.  We note also that a
candidate for a left-primitive Banach algebra which is not
right-primitive appeared in \cite{BoDu1}.  The Banach
algebras $A$ and $B$ considered below are much more tractable.

\section{Some representation theory for the Banach algebras $A$ and $B$}
\label{sec rep thy B alg}
\noindent
As in the introduction, $A$ is the Banach algebra $\ell^1(S)$ which
we may regard as all elements of the form $a=\sum_{m=0}^\infty x^m
\phi_m$ with $\phi_m\in\ell^1(\mathbb{Z})$, where we identify $\ell^1(\mathbb{Z})$ with
the Wiener algebra $W$.  We have $\|a\|=\sum_{m=0}^\infty \|\phi_m\|_1$
the latter norm being the $\ell^1$ norm.  Also, $B$ is the Banach
algebra of all elements of the form $b=\sum_{m=0}^\infty x^m \phi_m$
with $\phi_m\in C(\mathbb{T})$ and with the norm $\|b\|=\sum_{m=0}^\infty \|\phi_m\|_\infty$.

Our first task is to show that $A$ and $B$ are left-primitive. It is
not difficult to give continuous extensions of Irving's faithful
irreducible left representation of $\mathbb{C}S$ to both $A$ and $B$.
In passing we are able to simplify dramatically the provision in \cite{BoDu1} of an
example of specific dual representation behaviour.
We then consider an
averaging construction on irreducible right representations and
thereby obtain equivalent irreducible right representations on
classical Banach spaces. This enables us to describe all the
irreducible matrix right representations of $A$ and $B$.
Also, we introduce a family of non-faithful irreducible
infinite-dimensional right modules parameterized by a rich collection
of subsets of $\mathbb{T}$, namely infinite compact
sets which are square-closed, and minimal with respect to set
inclusion.  The argument for the case of $A$ requires us to
generalize the classical Wiener Lemma on the invertibility of
continuous functions on the circle with absolutely convergent
Fourier series.

We recall a fundamental construction for dual
representations of Banach algebras (see \cite{BoDu1} or \cite{BoDu2}).  Let
$\mathcal{A}$ be any Banach
algebra; for convenience, we suppose $\mathcal{A}$ has a unit, $1$.
Given $f\in\mathcal{A}'$ we get left and right ideals given by
$$
L_f = \{a: f(\mathcal{A}a) = (0)\} \quad \mbox{and} \quad K_f = \{a:f(a\mathcal{A}) = (0)\}.
$$
We have associated left and right
regular representations on the quotient spaces
${X_f=\mathcal{A}/L_f}$ and ${Y_f = \mathcal{A}/K_f}$, respectively.
The quotient spaces $X_f$ and $Y_f$ are in normed duality with\break
${\langle a',b' \rangle_f = f(ba)}$ for any choice of $a$ and $b$ in the cosets, and the
representations are linked by
$$
\langle a\xi,\eta\rangle_f = \langle\xi,\eta a\rangle_f.
$$
The left and right representations have the same kernel
$\{a:f(\mathcal{A}a\mathcal{A}) = (0)\}$.  The representations are
irreducible if and only if $L_f$ and $K_f$, respectively, are maximal one-sided ideals.
Every irreducible left representation of $\mathcal{A}$ is
equivalent to the representation on some $X_f$ (but then we know
nothing in general about the irreducibility of the right
representation on $Y_f$), and similarly for right representations.

We write $\delta_n$ for the usual point mass function on
$\mathbb{Z}$ (with value 1 at the point $n$ of
$\mathbb{Z}$). Irving's faithful irreducible left
representation of $\mathbb{C}S$ on $\mathbb{C}\mathbb{Z}$
is then determined by
$$
x \SP \delta_{2n-1} = 0, \quad x \SP \delta_{2n} = \delta_n, \quad y \SP \delta_n =\delta_{n-1}.
$$
It is routine to verify that the above
formulae determine a bounded left representation of $A$
on $\ell^1(\mathbb{Z})$. We write $V_A$ for
$\ell^1(\mathbb{Z})$ with this left module action.    For
the case of $B$ we rewrite $\delta_n$ as the function
$\zeta^n$ on the unit circle.  We take the corresponding left module
for $B$ to be $C(\mathbb{T})$.  The left action by $y$ is given
by
$$
y \SP f(\zeta) = \zeta^{-1}f(\zeta) \quad (\zeta\in\mathbb{T}).
$$
More generally the action by
$\phi(y)$ is given by
$$
\phi(y) \SP f(\zeta) = \phi(\zeta^{-1})f(\zeta) \quad (\zeta\in\mathbb{T}).
$$
The left action
by $x$ is well-defined by the formula
$$
x \SP f(\zeta) = \frac{1}{2} [f(\zeta^{1/2}) + f(-\zeta^{1/2})] \quad (\zeta\in\mathbb{T}).
$$
We denote this left module by $V_B$.

Here we introduce the natural dual
representations (see \cite[page 141]{BoDu2}) associated with $V_A$ and $V_B$.
There is a dual pairing on $V_A\times V_A$ given by
$$
\big\langle\sum \alpha_n \delta_n, \sum \beta_n \delta_n\big\rangle = \sum \alpha_n\beta_n
$$
and there is a dual
pairing on $V_B\times V_B$ given by
$$
\langle f,g\rangle = \int_0^{2\pi}f(e^{i\theta}) g(e^{-i\theta}) d\theta.
$$
In each case we get a dual
representation in which the right actions are (both) given by
\begin{equation}\label{eqn f zeta x phi}
f(\zeta)x = f(\zeta^2), \quad f(\zeta)\phi(y) = \phi(\zeta)f(\zeta).
\end{equation}
In the case of $A$,
(\ref{eqn f zeta x phi}) holds for all $f\in W$,
so we may identify $V_A$ with $W$.
For each module we have
$$
\langle a \SP v, w\rangle = \langle v,w \SP a\rangle
$$
and so the left representation is faithful if and only if the
right representation is faithful.  This allows us to choose
whichever leads to a simpler argument.

\begin{THM}\label{thm A B l-prim}
$A$ and $B$ are left-primitive Banach algebras.
\end{THM}

\begin{PRF}
For the left module $V_A$ we again write
$\{\delta_j\}$ for the usual normalized basis. Since
$y^k\delta_0=\delta_{-k}$ $(k\in\mathbb{Z})$
it follows that $\delta_0$ is strictly cyclic. We
show next that any non-zero $\xi=\sum\xi_j\delta_j$ can be mapped arbitrarily close to
$\delta_0$ by the left module action. It then follows that $A$ is
topologically irreducible on $V_A$. Since it has one strictly cyclic
vector it is then strictly irreducible on $V_A$ (see, for example,
\cite[Lemma~1]{CrMu}). We have $\xi_k\neq 0$ for some $k \in \mathbb{Z}$.
Note that $y^k \xi = \sum \xi_n \delta_{n-k}$ and that the
coefficient of $\delta_0$ is $\xi_k$, which is non-zero. Since
$\xi\in\ell^1(\mathbb{Z})$, it is straightforward to show that
$\xi_k^{-1}x^N y^k \xi \to \delta_0$ as $N\to\infty$.

To prove faithfulness we shall use the right representation. Thus we have
$$
\delta_n \SP x = \delta_{2n}, \quad \delta_n \SP y = \delta_{n+1}.
$$
Suppose that $V_A \SP a=0$ with $a = \sum_{m=0}^\infty x^m \phi_m(y)$. This
is equivalent to saying that $\delta_j \SP a=0$ for all
$j\in\mathbb{Z}$, that is, $\sum_{m=0}^\infty \delta_{2^m j}
 \SP \phi_m(y)=0$. Rewrite this in the notation of the Wiener algebra
$W$ and we have $\sum_{m=0}^\infty\zeta^{2^m j}\phi_m(\zeta) = 0$ for
all $\zeta\in\mathbb{T}$ and all $j\in\mathbb{Z}$.
Alternatively we
may write this as
$$
\sum_{m=0}^N \zeta^{2^m j}\phi_m(\zeta) + M(\zeta) = 0
$$
where $\|M\|_1$ is arbitrarily small. Multiply through by $\zeta^{-j-k}$,
integrate around the circle and let $|j| \to \infty$ to get
the $k$-th Fourier coefficient of $\phi_0$ arbitrarily small and
hence zero for all $k$.  This forces $\phi_0=0$.  Now repeat the
argument with a slight modification to get $\phi_1=0$. Similarly we
get all $\phi_m=0$ and so $a=0$, as required.

The above arguments require only minor modifications for the
module $V_B$.
\end{PRF}

\begin{RMK}
Let $T$ be the semigroup generated by $x$ and $y$
subject only to the relation $yx=xy^2$; thus we may regard $T$ as a
subsemigroup of $S$. Irving \cite{Ir1} proved that $\mathbb{C}T$ is both
left-primitive and right-primitive. He also proved that all
non-faithful (left and right) irreducible representations of
$\mathbb{C}T$ are finite-dimensional. We get Banach algebra
completions of $\mathbb{C}T$ by replacing $A$ by $\ell^1(T)$ and
modifying $B$ by replacing $C(\mathbb{T})$ by the disc algebra
$\mathcal{A}(\mathbb{D})$.  It may be verified that each completion
is both left-primitive and right-primitive; the modules are then
$\ell^1(\mathbb{Z}_+)$ and $\mathcal{A}(\mathbb{D})$.
\end{RMK}

In \cite{BoDu1} the authors presented a rather complicated example of a
dual representation of a Banach algebra in which the left
representation is irreducible while the right representation is
invariant on a chain of closed subspaces with zero intersection.
We can now present a simpler example of this behaviour.

\begin{THM}\label{thm dual rep A B}
For the dual representation of both $A$
and $B$ as given above, the left representation is
irreducible and the right representation is invariant on a chain of
closed subspaces with zero intersection.
\end{THM}

\begin{PRF}The argument is essentially the same in each case. We
present it for $A$.  For $p\in\mathbb{N}$, define
$$
V_p = \{f\in V_A:\mbox{$f(\zeta)=0$ whenever $\zeta^{2^p}=1$}\}.
$$
Clearly
each $V_p$ is a closed subspace and is non-zero since it contains
polynomials.  Also $V_{p+1} \suBset V_p$ and
$\cap_{p\in\mathbb{N}}V_p = (0)$ since any $f$ in the intersection
vanishes on a dense subset of $\mathbb{T}$.  Clearly $V_p$ is
invariant under right action by $y$ and it is also invariant under
$x$ since $\zeta^{2^{p+1}}=1$ whenever $\zeta^{2^p}=1$.
Then $V_p$ is invariant under right action by $A$, as required.
\end{PRF}

We turn now to an averaging construction for irreducible right representations.
Recall that any irreducible module for a Banach algebra
may be assumed to be a Banach space,
with continuous action.
Let the Banach space $V$ be an irreducible right
module for $A$ or $B$.  We may suppose without loss that $y$ and
$y^{-1}$ have operator norm $1$ on $V$ (apply \cite[Theorem 4.1]{BoDu2}), and
hence the operator $y$ has spectrum contained in $\mathbb{T}$.
{\it Suppose\/} that $v \SP y = \alpha v$ for some $\alpha\in\mathbb{T}$ and
some non-zero $v\in V$.
Choose $\rho\in V'$ with $\rho(v)=1$, and define
$f\in A'$ (or $f \in B'$) by
$$
f(a) = \LIM_{n\to\infty} \alpha^{-n}\rho(v \SP ay^n)
$$
where $\LIM$ is a Banach limit.
Note that
$f(1)=1$ and, for all $a$ in $A$ or $B$,
$$
f(ay)=\LIM_{n\to\infty} \SP \alpha^{-n}\rho(v \SP ayy^n)= \LIM_{n\to\infty} \alpha\alpha^{-n-1}\rho(v \SP ay^{n+1}) = \alpha f(a).
$$
We easily verify that
$$
K_f = \{a: f(ax^my^n)=0 \SP (m\in\mathbb{Z}_+,\ n\in\mathbb{Z})\}
=\{a: f(ax^m)=0 \SP (m\in\mathbb{Z}_+)\}.
$$
Let $K = \{a: v \SP a=0\}$ so that $K$ is a
maximal right ideal giving a right regular representation equivalent
to the given representation.  Clearly $f(K)=(0)$ and so $K\suBset
K_f$.  But $K_f$ is proper since $f(1)=1$, and hence $K=K_f$. For
each $a$ in $A$ (or $B$), define the bounded sequence on
$\mathbb{Z}_+$ by $g_a(m) = f(ax^m)$, and let $G$ be the space of
all such $g_a$. It is straightforward to check that our given right representation
is equivalent to the right representation on $G$ defined by $g_a \SP b =
g_{ab}$ (in fact, we simply map $V$ to $G$ by $v \SP a \to g_a$).
Moreover we have explicit formulae for the actions of $x$ and $y$ on $G$:
\begin{equation}\label{eqn g sub a y}
\left.
\begin{array}{c}
(g_a \SP x)(m) = g_{ax}(m) = f(axx^m) = g_a(m+1),
\\ [2ex]
(g_a \SP y)(m) = g_{ay}(m) = f(ayx^m) = f(ax^m y^{2^m}) = \alpha^{2^m}g_a(m).
\end{array}
\quad\right\}
\end{equation}
Notice what happens for the case  $\alpha=1$. Then
$y$ acts as the identity operator on $G$ and so the irreducible
image Banach algebra is commutative; but its centre is
$\mathbb{C}$ and hence $V$ is $1$-dimensional.  In general,
the $1$-dimensional representations come from the multiplicative
linear functionals $\chi$, and we pause here to list them.  They are
determined by the complex numbers $\chi(x)$ and $\chi(y)$ in the
closed unit disc subject to
$$
\chi(y)\chi(x) = \chi(yx) = \chi(xy^2) = \chi(x)\chi(y)^2.
$$
Since $y$ is invertible, we have either
$\chi(y)=1$ and $0<|\chi(x)|\le 1$, or $\chi(x)=0$ and $|\chi(y)|=1$.

We return to the general case of $\alpha\in\mathbb{T}$ and show that
this leads to the finite-dimensionality of the module $V$.  We have
already dealt with the $1$-dimensional case.

\begin{LEM}\label{lem irred r mod A B eigenvalue}
Let $V$ be an irreducible right
module for $A$ or $B$ with $\dim V\ge 2$, and
let $\alpha$ be an eigenvalue for $y$ on $V$ (and so
$\alpha\in\mathbb{T}$). Then $ \alpha^{2^k} = \alpha$
for some $k\in\mathbb{N}$ and $\dim V < \infty$.
Moreover, $x^k$ acts as a non-zero multiple of
the identity on $V$.
\end{LEM}

\begin{PRF}
We use the equivalent representation on $G$ which we
constructed above---see (\ref{eqn g sub a y}).
We have a non-zero $g$ with $g \SP y=\alpha g$ and hence
$$
\alpha g(m) = \alpha^{2^m} g(m) \quad (m\in\mathbb{Z}_+).
$$
If, for all $m\geq 1$, $\alpha^{2^m}\neq\alpha$ then we have $g(m)=0$ for
all $m\ge 1$.
Then $G=\mathbb{C}g$ and the representation is $1$-dimensional.
So we must have $\alpha^{2^k}=\alpha$ for some $k\in\mathbb{N}$.
It is routine to verify that $yx^k=x^k y$ as operators on $G$. It follows
that $x^k$ is in the centre of an irreducible Banach algebra of
operators and hence is a complex multiple of the identity.
Since $y^{2^k}=y$ as operators on $G$,
it follows that the image algebra is a finite-dimensional
irreducible algebra and hence $V$ is also finite-dimensional.
Let $k$ be minimal with the property that $\alpha^{2^k}=\alpha$.
Then $g(m)\neq 0$ implies $\alpha^{2^m}=\alpha$ which implies that $k$ divides $m$.
Suppose that $x^k=0$ on $G$.
For $m\geq k$, $g(m)=g \SP x^k(m-k) = 0$.
Hence $g(m)\neq 0$ implies $m=0$ and again the representation is $1$-dimensional.
\end{PRF}

Note that if $V$ is any finite-dimensional irreducible
right module for $A$ or $B$ then $y$ automatically has an eigenvalue
on $V$.  This enables us to list all the finite-dimensional
irreducible right representations for $A$ and $B$.  We determine all
the corresponding canonical matrices for $x$ and $y$.  These
matrices essentially appear in \cite{Ir1}.

We may suppose that $k$ is minimal with the property that
$\alpha^{2^k}=\alpha$.   It follows that the numbers
$\alpha,\alpha^2,\alpha^4,\ldots,\alpha^{2^{k-1}}$ are then
distinct. For, otherwise, $k \geq 2$ and we get $\alpha^{2^i}=\alpha^{2^j}$ with
$1\le i<j\le k$.  Raise both sides to the power $2^{k-j}$ and we
find that $\alpha^{2^{k+i-j}}=\alpha$ with $1\le k+i-j<k$, which is
impossible. We continue to work with the equivalent representation
on the subspace $G$ of $\ell^\infty$, so that $x^k = \beta 1$ on $G$
for some $|\beta|\le 1$. Consider first the case $\beta =1$. Since
$x^k=1$ on $G$ it follows that each $g\in G$ is periodic with period $k$.
Let $K=\{\alpha^{2^m}: m\in\mathbb{Z}_+ \}$, a finite set.
For $g\in G$, let $h(\alpha^{2^m})=g(m)$, and let $H$ denote all such
functions $h$. This is well-defined because of the periodicity of
each $g$. It is easy to verify that we get an equivalent right
representation on $H$ defined by
$$
(h \SP x)(\alpha^{2^m})=(g \SP x)(m)=g(m+1)=h(\alpha^{2^{m+1}}),
$$
$$
(h \SP y)(\alpha^{2^m})=(g \SP y)(m)=\alpha^{2^m}g(m)=\alpha^{2^m}h(\alpha^{2^m}).
$$
When we write these in standard function terms we get, for $\zeta\in K$,
\begin{equation}\label{eqn f zeta x zeta}
f(\zeta) \SP x = f(\zeta^2), \quad f(\zeta) \SP y = \zeta f(\zeta)
\end{equation}
---{\it cf}. (\ref{eqn f zeta x phi}).
By Theorem~\ref{thm M-set rep A B irred}, below, this right representation is irreducible for any
such $K$.
Alternatively, a computational proof of irreducibility can be constructed using the matrices below.
First verify that
$H$ contains the characteristic functions of the
singletons $\{\alpha^{2^j}\}$ for $j=0,1,\ldots, k-1$.
Take these functions as a basis for $H$, and
let $X$ and $Y$ be the matrices (with right action)
corresponding to the right actions of $x$ and $y$, respectively, on $H$.
It is straightforward to check that
$$
X= \left[ \begin{array}{ccccc}
0 & 0 & \cdots & 0 & 1 \\
1 & 0 & \cdots & 0 & 0 \\
0 & 1 & \cdots & 0 & 0\\
%\multicolumn{5}{c}\dotfill \\
\vdots & \vdots & & \vdots & \vdots \\
0 & 0 & \cdots & 1 & 0
\end{array} \right] ,
\quad Y = \left[ \begin{array}{ccccc}
\alpha & 0 & 0 & \cdots & 0 \\
0 & \alpha^2 & 0 & \cdots & 0 \\
0 & 0 & \alpha^4 & \cdots & 0 \\
\vdots & \vdots & \vdots & & \vdots \\
0 & 0 & 0 & \cdots & \alpha^{2^{k-1}}
\end{array} \right] .
$$

Return now to the general case with $x^k=\beta 1$.  Write
$\beta=\gamma^k$.  The right representations with $x$ and $y$ mapped to
$\gamma X$ and $Y$, respectively, are equally all irreducible with $0<|\gamma|\le 1$.
Write $\pi_{\alpha,\gamma}$ for this right representation.

\begin{THM}\label{thm fin dim r rep A B}
Up to equivalence, the irreducible
finite-dimensional right representations of $A$ and $B$
are given by the $1$-dimensional representations listed above
and the family $\pi_{\alpha,\gamma}$ (subject to
$\alpha^{2^k}=\alpha$ for some $k\in\mathbb{N}$, and
$0<|\gamma|\le 1$).  These representations separate the points
of $A$ and $B$.
\end{THM}

\begin{PRF}
It remains to prove the final assertion.  Suppose that
$a \SP \pi_{\alpha,\gamma} = 0$ for all these representations.  Let
$a=\sum x^m \phi_m(y)$.  Then we have
$$
\sum \gamma^m X^m \phi_m(Y) = 0
$$
for all complex $\gamma$ with $0<|\gamma|\le 1$.  It follows
that $X^m \phi_m(Y)=0$ for each $m$.  But $X$ is invertible and so
$\phi_m(Y)=0$ for each $m$.  We now have
$\phi_m(\alpha)=0$.
Since $\alpha$ can be any $(2^k-1)$-th root of unity,
this shows that $\phi_m=0$, giving $a=0$, as required.
\end{PRF}

\begin{COR}
Both $A$ and $B$ have the
property of direct finiteness.
\end{COR}

\begin{PRF}
Let $ba=1$ in $A$ or $B$.  Let $\pi$ be any
finite-dimensional irreducible right representation from the above
family.  Then $\pi(b)\pi(a)=\pi(1)=I.$  Since $\pi(a), \pi(b)$ are
complex matrices, we have $\pi(a)\pi(b)=I$.  Thus $\pi(ab-1)=0$ for
all such $\pi$, and hence $ab=1$ as required.
\end{PRF}

It is clear by dual representation theory that the family of
irreducible finite-dimensional left representations of $A$ and $B$
may be parameterized by mapping $x$ and $y$ to the matrices $\gamma X$ and $Y$,
respectively, where the matrices now act on the left on column vectors.

Suppose now that $K$ is any compact subset of $\mathbb{T}$ which is
square-closed, that is, $\zeta^2\in K$ whenever $\zeta\in K$.
It is easily verified that
$$
f(\zeta) \SP x = f(\zeta^2), \quad f(\zeta) \SP y = \zeta f(\zeta)
$$
give a right representation of $B$ on $C(K)$---{\it cf}. (\ref{eqn f zeta x zeta}).

We can think of $C(K)$, above, being derived from $C(\mathbb{T})$ by
the restriction mapping ${f \to f|_K}$.  But it is better to
regard it as the quotient of $C(\mathbb{T})$ by the submodule
$N=\{f\in C(\mathbb{T}): f(K)=(0) \}$.  To get the corresponding
right module for the Banach algebra $A$, we start with the analogous
right module $W$ and take the submodule $N=\{f\in W: f(K)=(0) \}$.
We then get the quotient module $W/N$ which we denote by $W(K)$.  Of
course $W(K)$ can be identified with the restrictions to $K$ of all
the functions in $W$.
Since $N$ is a closed ideal in the Banach algebra $W$
we can also consider $W(K)$ as a Banach algebra. The classic Wiener
Lemma adapts to $W(K)$.

\begin{LEM}\label{lem f in W(K) inv iff}
Let $f\in W(K).$
Then $f$ is invertible in $W(K)$ if and only if
$f$ is never zero on $K$.
\end{LEM}

\begin{PRF}
This can be deduced from results in \cite[Section~4.1]{Da}.
Alternatively, the argument in \cite[Example 19.4]{BoDu2} adapts as follows.
Write $u(\zeta) = \zeta$ for $\zeta\in\mathbb{T}$, and write $u|_K$ for
the restriction of $u$ to $K$.
It is enough to show that
$\SPEC(W(K),u|_K)$, the spectrum of $u|_K$ in $W(K)$, is $K$.
It is clear that $K\suBset \SPEC(W(K),u|_K)
\suBset \mathbb{T}$.
For each $\beta\in\mathbb{T}\setminus K$, we
have to show that $u|_K-\beta$ is invertible in $W(K)$.  Let the
interval $J$ be the component of the open set $\mathbb{T}\setminus
K$ which contains $\beta$.  We can modify the function $u-\beta$
only on $J$ so that it never vanishes on $J$, and we can do so
infinitely smoothly.  This modified function then has an absolutely
convergent Fourier Series and never vanishes on $\mathbb{T}$, and so
it has an inverse in $W$ by Wiener's Lemma.  Restrict this inverse
to $K$ and we see that $u|_K-\beta$ is invertible in $W(K)$.
\end{PRF}

When $K$ is a proper subset of $\mathbb{T}$, the kernel of this
right representation of $B$ intersects the subalgebra
$C(\mathbb{T})$ of $B$ in the space of all functions that vanish on
$K$, and so this representation is certainly not faithful. A similar
remark applies for $A$. But we show below that the representation is
irreducible if and only if $K$ is minimal, with respect to set
inclusion, amongst all compact square-closed subsets of
$\mathbb{T}$.

\begin{THM}\label{thm M-set rep A B irred}
Let $K$ be a compact
square-closed subset of $\mathbb{T}$. Then the above right
representation of $B$  (respectively $A$) is irreducible on
$C(K)$ (respectively $W(K)$) if and only if $K$ is minimal.
\end{THM}

\begin{PRF}
If $K$ is not minimal, then any smaller compact square-closed set
produces an invariant subspace for the representation.  Suppose now
that $K$ is minimal. The function 1 is cyclic.
Let $f$ be any non-zero function in $C(K)$ and
let $g=|f|^2= f f^*$
(here, $f^*(\zeta)$ is the complex conjugate of $f(\zeta)$).
Then $g\in fB$. Let $h = \sum_{n=0}^\infty 2^{-n} g \SP x^n$.
Then $h\in fB$. Note that
$$
h(\zeta) = \sum_{n=0}^\infty 2^{-n} g(\zeta^{2^n})
$$
so that $h$ is non-negative and not identically zero.
Let $L = \{\zeta\in K:h(\zeta)=0 \}$ so that $L$ is a compact proper subset of $K$.
Let $\zeta\in L$. Then $g(\zeta^{2^n})=0$ for each $n$ and so
$h(\zeta^2)=0$, that is, $\zeta^2\in L$. Since $K$ is minimal, it
follows that $L$ is empty and so $h$ is invertible in $C(K)$.  Thus
$1\in hB \suBset fB$ and so we can map $f$ to 1.  The proof is
complete for the case of $B$. The proof adapts for the case of $A$,
with the use at the last step of the above Wiener Lemma for $W(K)$.
\end{PRF}

For brevity we shall refer to minimal compact
square-closed subsets of $\mathbb{T}$ as $\mu$-{\it sets}.  The
existence of infinite $\mu$-sets is guaranteed by Zorn's Lemma, but that
approach gives us little insight.
We have already met a countable family of finite $\mu$-sets:
the sets of the form $\{\alpha^{2^m}:m\in\mathbb{Z}_+\}$
where $\alpha^{2^k}=\alpha$ for some $k\in\mathbb{N}$.
There are in fact uncountably
many different infinite $\mu$-sets (and each one is a Cantor set). Two
such uncountable families are given in \cite{Ve}.
The union of all infinite $\mu$-sets is dense in $\mathbb{T}$, and so the
corresponding family of irreducible infinite-dimensional right
representations again separates the points of $A$ and $B$.  This
density property of $\mu$-sets does not appear to be in the literature,
but the authors give a proof in~\cite{CrDuMc}.

We conjecture that the only irreducible right representations
of dimension greater than $1$ for
$B$ are those given by minimal $K$ as above.  We do not know whether
there exist any non-faithful irreducible left representations of $A$
or $B$ on infinite dimensional spaces.

\section{Appendix: some representation theory for the semigroup algebra $\mathbb{C}S$}
\label{sec appendix}
\noindent
We present here a different proof of Irving's result in \cite{Ir2} that
$\mathbb{C}S$ is left-primitive but not right-primitive.
The proof of Lemma~\ref {lem irred r mod CS eigenvalue}
is largely a translation
of Irving's arguments with ideals, to arguments with vectors.
We prove further that the only irreducible right modules of $\mathbb{C}S$ are
finite-dimensional.  As before $S$ is the semigroup with identity, generated by
$x$ and $y$, subject to $y$ being invertible and $yx=xy^2$.  Thus the
product in $S$ is given by
$$
x^my^n \SP x^py^q = x^{m+p}y^{n2^p+q}.
$$
There is a natural left module $V_S$ for $\mathbb{C}S$ given by
taking $V_S$ to be the linear space $\mathbb{C}\mathbb{Z}$ and
restricting to $V_S$ the module action for $V_A$ defined in Section~\ref{sec rep thy B alg}.
To see that $\mathbb{C}S$ is left-primitive, we simply adapt the
proof of Theorem~\ref{thm A B l-prim}.
The use of Fourier analysis in
proving faithfulness can be avoided---we leave the reader to give a direct
algebraic proof that $a \SP V_S\neq(0)$ whenever $a\neq 0$.
We record the result.

\begin{THM}\label{thm CS l-prim}
The linear algebra $\mathbb{C}S$
is left-primitive.
\end{THM}

The representation theory for $\mathbb{C}S$ is of course different
from the representation theory for $A$ since there are no continuity
restrictions for $\mathbb{C}S$.  For example, there are two (larger)
families of $1$-dimensional irreducible representations of
$\mathbb{C}S$ determined by multiplicative linear functionals
$\chi$.  As before, these are determined by $\chi(x)$ and $\chi(y)$
subject to
$$
\chi(y)\chi(x)=\chi(x)\chi(y)^2.
$$
Since there is no
boundedness constraint on $\chi$, we have either $\chi(y)=1$ and
$\chi(x)\in\mathbb{C}$, or $\chi(x)=0$ and $\chi(y)\in\mathbb{C}$
with $\chi(y)\neq 0$.  Similarly we get a larger two-parameter family
of finite-dimensional irreducible right representations of
$\mathbb{C}S$.  As before, let $\alpha\in\mathbb{C}$ with
$\alpha^{2^k} = \alpha$ and the positive integer $k$ minimal.  Now
let $\gamma\in\mathbb{C}$, $\gamma\neq 0$, and let $X$ and $Y$ be the
matrices as in Section~\ref{sec rep thy B alg}.
Let $\pi_{\alpha,\gamma}$ be the right
representation of $\mathbb{C}S$ determined by
$$
x \SP \pi_{\alpha,\gamma} = \gamma X, \quad y \SP \pi_{\alpha,\gamma} = Y.
$$
As for $A$ and $B$ in Section~\ref{sec rep thy B alg},
we now get the following result for $\mathbb{C}S$.

\begin{THM}\label{thm fin dim r rep irred}
The above finite-dimensional right
representations $\pi_{\alpha,\gamma}$ of $\mathbb{C}S$
are irreducible and the family separates the points of
$\mathbb{C}S$.
\end{THM}

\begin{COR}
$\mathbb{C}S$ has the property of direct
finiteness.
\end{COR}

\begin{RMK}
If we regard the matrices $X$ and $Y$ as acting on the left,
then we get a corresponding family of irreducible
finite-dimensional left representations of $\mathbb{C}S$ with the
corresponding properties.
\end{RMK}

To help simplify the notation in places, we shall write $R$
for the algebra of Laurent polynomials in $y$; thus $R$ is simply a
copy of $\mathbb{C}\mathbb{Z}$.  We aim to show that 
any irreducible right module $V$ for $\mathbb{C}S$ is
finite-dimensional, and it will follow immediately that
$\mathbb{C}S$ is not right-primitive.  The first (and most critical)
step is to show that $y$ must have an eigenvalue on V.  We can then
adapt the right averaging construction to our algebraic situation to
prove that $V$ must be finite-dimensional.  In the non-Banach
algebra setting we need a proof that an irreducible
representation of $\mathbb{C}\mathbb{Z}$ is $1$-dimensional.
We could use the {\it Nullstellensatz},
but in the course of Lemma~\ref{lem irred r mod CS eigenvalue}
below we give instead a very short elementary argument for this.

\begin{LEM}\label{lem irred r mod CS eigenvalue}
Let $V$ be an irreducible right
module for $\mathbb{C}S$.
Then $y$ has an eigenvalue on~$V$.
\end{LEM}

\begin{PRF}
Suppose first that $V \SQ x = (0)$. Let $v\in V$, $v\neq 0$.
Consider $v \SP (1+y)$. If $v \SP (1+y)=0$ then $y$ has an
eigenvalue.  If $v \SP (1+y)\neq 0$, then by irreducibility there is
a polynomial $\phi$ in $y$ and $y^{-1}$ such that $v \SP
(1+y)\phi=v$. Multiply by a suitable power of $y$ to get $v \SP
(1+y)\psi = v \SP y^k$ for some polynomial $\psi$ in $y$. Since
$y^k-(1+y)\psi$ is a non-zero polynomial it follows, by factorizing
into linear factors, that $y$ has an eigenvalue.

Now suppose that $V \SQ x \neq (0)$. Let $v\in V$ with $v \SP x \neq
0$. Since $V$ is irreducible, there is $d\in\mathbb{C}S$ with $v \SP
xd=v$. This gives
$$
0 = v \SP (1-xd) = v \SP (1+xd_1+x^2d_2+\cdots)
$$
for some $d_j\in R$. Suppose that $a=r_0+xr_1+\cdots+x^mr_m$ has $m$ minimal
such that $m\ge 0$, $va=0$, $r_j\in R$ and $r_0\neq 0$. Our aim is to
show that $m=0$. As above, this will prove that $y$ has an
eigenvalue on $V$.

From now on, we shall also regard the elements of $R$ as functions
on $\mathbb{T}$, via $y(\zeta)=\zeta$. We first choose
an odd prime $p$ so that for each zero $e^{2\pi i\beta_j}$ of $r_m$
with $0<\beta_j<1$ we have $p\beta_j \notin\mathbb{Z}$. For this we
need consider only those $\beta_j$ which are rational, and we may
take the prime $p$ to be greater than any odd prime which appears in
a denominator of any $\beta_j$ in its lowest terms.  Put $\alpha =
e^{2\pi i/p}$.  Then $\alpha^{2^k} \neq e^{2\pi i\beta_j}$ for any
$j$ or for any $k\in\mathbb{Z}_+$; for otherwise we should
have $2^k 2\pi /p - 2\pi\beta_j \in 2\pi\mathbb{Z}$ and hence
$p\beta_j\in\mathbb{Z}$. Also, $\alpha^{2^k}\neq 1$ since
$2^k/p\notin\mathbb{Z}$.   We thus have $r_m(\alpha^{2^k})\neq 0$ for
all $k\in\mathbb{Z}_+$.

Now put $q=y^p-1$ and let $Q=Rq$.  We have
$qx=x(y^{2p}-1)=xq(y^p+1)\in xQ$.  Similarly, $qx^2\in x^2Q$, $qx^3\in
x^3Q$, and so on.  If $v \SP q=0$ then $y$ has an eigenvalue on $V$ and we
are done. Suppose instead that $vq\neq 0$ and so there is
$d'\in\mathbb{C}S$ with $vqd'=v$. This gives $v(1-qd')=0$ where
$$
1-qd' = 1+q_0+xq_1+x^2q_2+\cdots
$$
and each $q_j\in Q$.
By the choice of $\alpha$ we have $q(\alpha)=0$ for every $q\in Q$.
Thus $(1+q_0)(\alpha)=1$.  Suppose now that
$$
b=t_0+xq_1+\cdots+x^nq_n
$$
has $n$ minimal such that $n\ge 0$,
$vb=0$, $q_j\in Q$ and $t_0\in R$ with $t_0(\alpha)\neq 0$.  Since $b$
satisfies the conditions given for $a$ above, it follows that $n\ge
m$.  Put $n=m+k$ so that $k\ge 0$.  We now suppose that $m>0$ and
derive a contradiction.  This will give $m=0$ and so complete the
proof of the lemma.  For $r\in R$ put $\hat{r}(y)=r(y^{2^k})$. With
$a$ as above, we have
$$
0=v \SP ax^k=v \SP (x^k \hat{r}_0+\cdots + x^n\hat{r}_m)
$$
and hence $0=v \SP c$ where
$$
c=b\hat{r}_m-ax^kq_n = T_0+xQ_1+\cdots + x^{n-1}Q_{n-1}
$$
with $T_0\in R$ and $Q_j\in Q$.  Here we have $T_0
= t_0\hat{r}_m$ if $k>0$, and $T_0=t_0 r_m-r_0 q_n$ if $k=0$. In
either case we have
$$
T_0(\alpha) = t_0(\alpha)\hat{r}_m(\alpha) = t_0(\alpha) r_m(\alpha^{2^k}) \neq 0.
$$
Thus $c$ has the same form as
$b$ but is of lesser degree, which is impossible.
Hence $m=0$, which was needed to complete the proof.
\end{PRF}

\begin{THM}\label{thm irred r mod CS fin dim}
Let $V$\BS be any irreducible right
module for $\mathbb{C}S.$
Then $V$\BS is
finite-dimensional, and the representation is $1$-dimensional (as
listed above) or is equivalent to one of the form
$\pi_{\alpha,\gamma}.$ Thus $\mathbb{C}S$ is not
right-primitive.
\end{THM}

\begin{PRF}
By Lemma~\ref {lem irred r mod CS eigenvalue}
we have $v\in V$ and $\alpha\in\mathbb{C}$ such that
$v \SP y = \alpha v$.
Since $y$ is invertible, $\alpha\neq 0$.
Let $\rho$ be a linear functional on $V$ with
$\rho(v)=1$.
Consider first the case $|\alpha|\geq 1$.
We note that $\alpha^{-n}\rho(v \SP y^n)=1 \SP
(n\in\mathbb{N})$.
Also
$$
\alpha^{-2m}\rho(v \SP xy^{2m}) = \alpha^{-2m}\rho(v \SP y^mx) = \alpha^{-m}\rho(v \SP x)
$$
and similarly $\alpha^{-2m-1}\rho(v \SP xy^{2m+1})=\alpha^{-m-1} \rho(vxy)$.
Thus $\{\alpha^{-n}\rho(v \SP xy^n): n\in\mathbb{N}\}$ is a bounded
sequence and hence we may apply a Banach limit to define
$$
f(x) = \LIM_{n\to\infty} \alpha^{-n} \rho(v \SP xy^n).
$$
Similarly we may define
$$
f(x^j) = \LIM_{n\to\infty} \alpha^{-n}\rho(v \SP x^jy^n) \quad (j\in\mathbb{N})
$$
and then we can extend to a
linear functional $f$ on $\mathbb{C}S$ given by
$$
f(a) = \LIM_{n\to\infty} \alpha^{-n} \rho(v \SP ay^n).
$$
If $|\alpha| < 1$ we instead define $f(a) = \LIM_{n\to\infty} \alpha^n \rho(v \SP ay^{-n})$.
In each case, as before,
$f(1)=1$ and $f(ay) = \alpha f(a)$.
We can now proceed exactly as in Section~\ref{sec rep thy B alg}
to get our given right representation
equivalent to the right representation on the space
$G=\{g_a:a\in\mathbb{C}S\}$ with the actions $(g_a x)(m)=g_a(m+1)$
and $(g_a y)(m)=\alpha^{2^m}g_a(m)$ for $m\in\mathbb{Z}_+$.

We have $h\in G$, $h\neq 0$, with $h \SP y = \alpha h$.
If, for all $m\in\mathbb{N}$, $\alpha^{2^m}\neq\alpha$ then $h(m)=0$ for all $m\in\mathbb{N}$.
Then $h \SP x = 0$.
In this case $G=\mathbb{C}h$ and the representation is $1$-dimensional.

Finally, consider the case $\alpha^{2^k}=\alpha$ for some $k\in\mathbb{N}$.
We may assume that $k$ is minimal.
For simplicity we shall now consider the case $k=2$; the general
case is just an elaboration of the argument.  Thus we have
$\alpha^4=\alpha\neq \alpha^2$.  We write $X$ and $Y$ for the operators on
$G$ corresponding to $x$ and $y$, respectively. It is immediate from the action of $y$
given above that the algebra generated by $Y$ consists of all
diagonal operators with entries of the form
$(\gamma,\delta,\gamma,\delta,\dots)$.
It is clear that our
eigenvector $h$ for $y$ has to be of the form
$h=(h_0,0,h_2,0,h_4,0,\ldots)$.
If $h X^2 = 0$ then $h=(h_0,0,0,\ldots)$,
and again our representation is $1$-dimensional.
Suppose that $h X^2 \neq 0$.
By irreducibility we get
$hX^2\sum_jX^jp_j(Y)=h$ and since the support of $hX^{2m+1}$ is
contained in $\{1,3,5,\ldots\}$, we can assume that only even $j$
appear in the sum.  But each $p_{2m}(Y)$ acts as a multiple of $I$
on $hX^{2+2m}$.  It follows that $hQ(X)=0$ for some non-zero complex
polynomial $Q$.
The usual argument shows that $X$ has
an eigenvector $g'\in G$. Suppose that $g' X=\beta g'$ for some
$\beta\in\mathbb{C}$.
Then we have $g'=(1,\beta,\beta^2,\beta^3,\ldots)$.
If $\beta = 0$ the representation is $1$-dimensional.
Suppose $\beta\neq 0$.
Let $G_1$ be the $2$-dimensional subspace of $G$ given by
$$
G_1 = \{ g\in G: g(j+2)=\beta^2 g(j)\ (j\in\mathbb{Z}_+)\}.
$$
A calculation shows that $G_1$ is invariant
under $X$ and $Y$. Since $G_1\neq (0)$, it follows by irreducibility that
$G_1=G$ and hence $\dim G=2$.

Verification that our representation is equivalent
to the representation $\pi_{\alpha,\beta}$ is routine.
\end{PRF}

\begin{RMK}
The use of $\LIM$ can be avoided in Section~\ref{sec rep thy B alg}
by using numerical range results from \cite{Si};
in Section~\ref{sec appendix} the $\LIM$s can be calculated explicitly, if desired.
\end{RMK}

\end{document}